\documentclass[letterpaper,10pt,conference]{ieeeconf}  %
\usepackage{latexsym,amssymb}

\IEEEoverridecommandlockouts

\newcommand{\comment}[1]{}
\newcommand{\R}{{\mathbb R}}  

\newtheorem{theorem}{Theorem}
\newtheorem{itlemma}{Lemma}[section] 
\newtheorem{itproposition}[itlemma]{Proposition}
\newtheorem{itcorollary}[itlemma]{Corollary}
\newtheorem{itremark}[itlemma]{Remark}
\newtheorem{itdefinition}[itlemma]{Definition}
\newtheorem{itexample}[itlemma]{Example}

\newenvironment{lemma}{\begin{itlemma}\rm}{\end{itlemma}} 
\newenvironment{remark}{\begin{itremark}\rm}{\end{itremark}}
\newenvironment{corollary}{\begin{itcorollary}\rm}{\end{itcorollary}}
\newenvironment{proposition}{\begin{itproposition}\rm}{\end{itproposition}}
\newenvironment{definition}{\begin{itdefinition}\rm}{\end{itdefinition}}
\newenvironment{example}{\begin{itexample}\rm}{\end{itexample}}

\newcommand{\text}[1]{\hbox{\rm \ #1\ \/}}
\newcommand{\be}[1]{\begin{equation}\label{#1}}
\newcommand{\ee}{\end{equation}}
\newcommand{\bl}[1]{\begin{lemma}\label{#1}}
\newcommand{\ble}[1]{\begin{lemmaex}\label{#1}}
\newcommand{\br}[1]{\begin{remark}\label{#1}}
\newcommand{\bt}[1]{\begin{theorem}\label{#1}}
\newcommand{\bd}[1]{\begin{definition}\label{#1}}
\newcommand{\bp}[1]{\begin{proposition}\label{#1}}
\newcommand{\bc}[1]{\begin{corollary}\label{#1}}
\newcommand{\bfact}[1]{\begin{fact}\label{#1}}
\newcommand{\ber}[1]{\begin{exercise}\label{#1}}
\newcommand{\bex}[1]{\begin{example}\label{#1}}
\newcommand{\bem}[1]{\begin{example}\label{#1}} 
\newcommand{\ec}{\mybox\end{corollary}}
\newcommand{\efact}{\mybox\end{fact}}
\newcommand{\eer}{\mybox\end{exercise}}
\newcommand{\eex}{\mybox\end{example}}
\newcommand{\eem}{\mybox\end{example}}
\newcommand{\el}{\mybox\end{lemma}}
\newcommand{\ele}{\mybox\end{lemmaex}}
\newcommand{\er}{\mybox\end{remark}}
\newcommand{\et}{\qed\end{theorem}}
\newcommand{\ed}{\mybox\end{definition}}
\newcommand{\ep}{\mybox\end{proposition}}
\newcommand{\epr}{\end{proof}}
\newcommand{\bpr}{\begin{proof}}

\newcommand{\ecs}{\end{corollary}}
\newcommand{\eers}{\end{exercise}}
\newcommand{\eexs}{\end{example}}
\newcommand{\eems}{\end{example}}
\newcommand{\els}{\end{lemma}}
\newcommand{\eles}{\end{lemmaex}}
\newcommand{\ers}{\end{remark}}
\newcommand{\ets}{\end{theorem}}
\newcommand{\eds}{\end{definition}}
\newcommand{\eps}{\end{proposition}}
\newcommand{\halmos}{\rule{1ex}{1.4ex}}
\newcommand{\qed}{\hfill \halmos} 
\newcommand{\mybox}{\hfill $\Box$} 
\newcommand{\beq}{\begin{eqnarray}}
\newcommand{\eeq}{\end{eqnarray}}
\newcommand{\beqn}{\begin{eqnarray*}}
\newcommand{\eeqn}{\end{eqnarray*}}
\newcommand{\bi}{\begin{itemize}}
\newcommand{\ei}{\end{itemize}}
\newcommand{\ben}{\begin{enumerate}}
\newcommand{\een}{\end{enumerate}}

\newcommand{\kk}{{\mathcal K}}
\newcommand{\kl}{{\cal KL}}
\newcommand{\ki}{{\cal K}_\infty}
\newcommand{\norm}[1]{\left\Vert#1\right\Vert}
\newcommand{\abs}[1]{\left\vert #1 \right\vert}
\newcommand{\rref}[1]{(\ref{#1})}
\newcommand{\edo}{\end{document}}

\newcommand{\iss}{ISS}

\newcommand{\xc}{{\cal X}}
\newcommand{\vc}{{\cal V}}
\newcommand{\vz}{{\cal V}_0}
\newcommand{\uc}{{\cal U}} 
\newcommand{\id}{{\bf id}}     
\newcommand{\veps}{{\varepsilon}}
\newcommand{\ve}{{\varepsilon}}

\title{\LARGE \bf A Nonlinear Small-Gain Theorem for Large-Scale Time
  Delay Systems}
\date{}

\author{Shanaz Tiwari, Yuan Wang, and Zhong-Ping Jiang
\thanks{This work was supported partially by NSF grants
 DMS-0504462, DMS-0504296 and the Chinese National 
Natural Science Foundation grants 60228003 and 60628302.
}
\thanks{S. Tiwari and Y. Wang are with the Department of Mathematical
  Sciences, Florida Atlantic University, Boca Raton, FL 33431.
        {\tt\small \{stiwari1, ywang\}@fau.edu}}%
\thanks{Z.P. Jiang is with the Department of Electrical and Computer
 Engineering,  Polytechnic Institute of New York University, 
 Brooklyn, NY 11201, USA.
{\tt\small zjiang@control.poly.edu}}%
}

\begin{document}

\maketitle
\thispagestyle{empty}
\pagestyle{empty}

\begin{abstract}
This paper extends the nonlinear ISS 
small-gain theorem to a large-scale time delay system composed of
three or more subsystems.  En route to proving this small-gain
theorem for systems of differential equations with delays, 
a small-gain theorem for operators
is examined.  The result developed for operators allows applications to
a wide class of systems, including state space systems with delays.

\end{abstract}

\section{Introduction}

One of the most powerful tools in stability analysis and control design
of interconnected systems is the small-gain theory.  The first
small-gain theorem in the context of input-to-state stability (ISS) was 
developed in \cite{jiang-praly-teel}.  A variant of the ISS
small-gain theorem was given in \cite{CPT} in terms of asymptotic gains.
The authors of \cite{is-smallgain} presented 
an \iss-type small-gain theorem in terms of operators 
that recovers the case of state space form, with applications to
several variations of the \iss\ property such as input-to-output
stability, incremental stability and detectability for interconnected
systems. Initialized by the work \cite{DRW07} and \cite{Teel03}, the
small-gain theory was extended to systems composed of three or more
subsystems. In the recent work \cite{JW-WCICA}, a cyclic small-gain
theorem was provided to deal with the input-to-output stability (IOS) 
properties for large-scale interconnected systems.

The goal of the current work is to develop the small-gain results for
large-scale systems with time-delays appearing in both the subsystems
and the interconnections.  Systems with delays arise naturally from many
practical applications such as networked control systems.  As a
consequence, time-delay systems and control have received much attention
in recent years, see for instance, \cite{fan-arck-2006},
\cite{jankovic-2001}, and \cite{MML-2008}.
In a series of recent work \cite{Kar07c},
\cite{Kar07d}, and \cite{Kar-08}, various notions related to ISS were
studied for systems with delays.  As in the case for systems without
delays, small-gain theorems provide natural tools for stability analysis
of interconnected systems.  In \cite{teel-delay}, a Razumikhin-type
theorem on stability analysis was presented for systems with delays by
using the nonlinear small-gain theorem.  In \cite{PTM-2006}, a
small-gain theorem was developed to solve a stabilization problem of a
force-reflecting telerobotics system with time delays. In
\cite{Kar-Jiang-07} a small-gain theorem was given for a wide class of
systems including systems with time delays.

The previous work on small-gain theorems for time-delay 
systems focused on systems composed of two subsystems.
Our main contribution will be to present a small-gain theorem for
interconnected systems composed of three or more subsystems
with time delays.  Instead of carrying out our proofs for systems of
differential equations with delays, we develop a small-gain
theorem for input/output operators.  This is an approach adopted in
\cite{is-smallgain}.  The advantage of doing so is that it allows one to
develop small-gain theorems for a wide class of systems including
systems of differential equations with delays and possibly certain types
of hybrid systems.

{\it Notations.} Throughout this work, we use $|\cdot|$ to denote the
Euclidean norm of vectors, and $\|\cdot\|_{I}$ to denote the
essential supremum norm of measurable and locally essentially bounded
functions defined on the interval $I$. 
For $\phi=(\phi_1,\cdots, \phi_k)$ defined on an interval $I$, we
let $\displaystyle{\norm{\phi}_{I}= \max_{1\leq i \leq
  k}\left\{\norm{\phi_i}_{I}\right\}}$. A function $\alpha :
\R_{\geq 0} \to \R_{\geq 0}$ is {\em of class $\kk$} if it is
continuous, positive definite, and strictly increasing; and is {\em of
  class $\ki$} if it is also unbounded.  A function $\beta : \R_{\geq 0}
\times \R_{\geq 0} \to \R_{\geq 0}$ is said to be of {\em class $\kl$}
if for each fixed $t \geq 0$, $\beta(\cdot, t)$ is of class $\kk$, and
for each fixed $s \geq 0$, $\beta(s,t)$ decreases to 0 as $t \to
\infty$.

\section{Preliminaries}

Let $\theta \ge 0$, and let $C[-\theta, 0]$ denote the Banach space of
continuous functions defined on $[-\theta, 0]$, equipped with the norm
$\norm{\,\cdot \,}_{[-\theta, 0]}$.  For a continuous function $q: [-\theta,
b)\rightarrow\R$, where $b > 0$, define $q_t(s):= q(t+s)$.  
Then, for each $t\in [0, b)$, $q_t\in C[-\theta, 0]$.

Consider a nonlinear system with time delays described by
\beq
(\Sigma_u):\ \dot{x}(t)&=& f\left(x_t,v_t,u(t)\right), \ \ \ t\geq 0,
\label{e-sys}
\eeq
where 
\begin{itemize}
\item for each $t$, $x(t)\in \R^n$, \ $v(t)\in \R^p$ and  $u(t)\in \R^m$; 
\item  $f:\xc \times \vz \times \R^m\rightarrow \R^n$ is locally
  Lipschitz and completely continuous (see \cite{Hale-book} for definition), 
  where $\xc = \left(C[-\theta,\ 0]\right)^n$ and 
$\vz = \left(C[-\theta,\ 0]\right)^p$.
\end{itemize}
The input of the system is given by $(v(t),\ u(t))$ where $v: 
[-\theta, \infty) \rightarrow \R^p$ 
is continuous and $u: [0, \infty)\rightarrow \R^m$
is measurable and locally essentially bounded. Trajectories of the
systems are absolutely continuous functions defined on some interval that
satisfy \rref{e-sys} almost everywhere.  We let $\uc$ denote the collection
of measurable and locally essentially bounded functions, and let $\vc =
  \left(C[-\theta, \infty)\right)^p$. 

With the assumptions on $f:\xc \times \vz \times \R^m\rightarrow \R^n$
given above, it can be shown that for each $v(\cdot) \in \vc$ and each
$u(\cdot) \in \uc$, the map $F: \xc \times \R_{\geq 0}\rightarrow \R^n$
given by 
\be{e-F}
F\left(\xi,\ t\right):= f\left(\xi,\ v_t,\ u(t)\right)
\ee
is bounded on any compact set of $\xc \times [0,\ \infty)$ and is
 locally Lipschitz in $\xi$, uniformly in $t$,
  for all $t$ in any compact
 set (see Lemma \ref{l-Lip} in Appendix A.) It can also be shown that 
 $F$ is completely continuous in $\xi$, uniformly in $t$,
  for all $t$ in any compact
 set. 
Consequently, for each $v\in\vc, u\in\uc$, and
each continuous function $\xi$ 
defined on $[-\theta, 0]$, there is a unique trajectory of \rref{e-sys}
corresponding to $v$ and $u$ that 
satisfies the initial condition $x_0(\cdot) = \xi(\cdot)$
(see \cite{Hale-book}).  We denote this 
 trajectory by $x(t,\xi,  v, u)$, and its maximum interval by
$[-\theta, \ T^{\max}_{\xi, v, u})$.

Our reason for considering $\dot{x}(t)= f\left(x_t, v_t, 
  u(t)\right)$ instead of the cases of $\dot{x}(t)=
f\left(x_t,u(t)\right)$ or $\dot{x}(t)= f\left(x_t, v_t\right)$ is 
that we want to allow an interconnected sytem
to have feedbacks involving time delays (c.f. Section \ref{s-3}), and
the input signals without delays to be merely measurable.  Also
note that the form of system \rref{e-sys} allows $f$ to depend on
  $x_t(\cdot)$ in any manner, which enables \rref{e-sys} to cover a wide
  class of time-delay systems, e.g., systems involving discrete
  time-delays as in 
$\dot x(t)
  = g(x(t), x(t-\theta_1), \ldots, x(t-\theta_l), v(t-\theta_0), u(t)),
$ ($0\le \theta_i\le \theta$)
or systems involving distributed time-delays as in
$$
\dot x(t) = G\left(\frac{1}{\theta}\int_{t-\theta}^t x(s)\,ds, \;
u(t)\right), \ \ (\theta > 0).
$$
 
 \bd{ISS} The system $(\Sigma_u)$ is (globally) {\it{input-to-state
 stable}} (ISS) if there exist $\kk$-functions $\gamma^u(\cdot)$
 and $\gamma^v(\cdot)$  and a $\kl$-function $\beta(\cdot)$ such that
 \beq 
& &\abs{x(t, \xi, v, u)} \leq  \beta \left(\norm{\xi}_{[-\theta,\ 0]}, t
 \right)\nonumber \\
 & &\quad  +\ \gamma^v\left(\norm{v}_{[-\theta,\infty)}\right) 
 + \, \gamma^u\left(\norm{u}_{[0,\infty)}\right)
\label{isseq}
 \eeq 
for all $t\ge 0$.
\eds

\bd{GS}The system $(\Sigma_u)$ is said to satisfy the {\it{global
 stability}} (GS) property 
 if there exist $\kk$-functions $\sigma^x(\cdot)$, $\sigma^v(\cdot)$ and
 $ \sigma^u(\cdot)$ such that for all $t\ge 0$,
\beq
\abs{x(t,\xi, v, u)} \!\!&\leq&\!\!
\max\left\{\sigma^x\left(\norm{\xi}_{[-\theta,\
 0]}\right),\right. \nonumber \\  
& &\left. \!\!\sigma^v\left(\norm{v}_{[-\theta,\infty)}\right), 
\sigma^u\left(\norm{u}_{[0,\infty)}\right)\right\}\!.
\label{gseq}
\eeq
\eds
 
\bd{AG} The system $(\Sigma_u)$ is said to satisfy the {\it{asymptotic
 gain}} (AG) property 
 if there are $\kk$-functions $\gamma^u(\cdot)$ and
 $\gamma^v(\cdot)$ such that 
\beq
& &\overline{\lim_{t\rightarrow \infty}}\abs{x(t,\xi, v, u)}\nonumber\\
& & \leq \max\left\{\gamma^u \left(\norm{u}_{[0,\infty)}\right), \ 
\gamma^v \left(\norm{v}_{[-\theta,\infty)}\right)\right\}\label{ageq}
\eeq
 for all $u \in \uc$ and  $v \in \vc$.
\eds

It can be seen that the \iss\ property implies both the (GS) and the (AG)
conditions. 

\bl{aglem}
The AG condition as given in \rref{ageq} is equivalent to 
\beq
& &
\overline{\lim_{t\rightarrow
    \infty}}\abs{x(t,\ \xi, \ v, \ u)}\nonumber\\
& &\leq\max\left\{\overline{\lim_{t\rightarrow
    \infty}}\gamma^u\left(\abs{u(t)}\right),  \
  \overline{\lim_{t\rightarrow
    \infty}}\gamma^v\left(\abs{v(t)}\right)\right\}.\label{eqivag}
\eeq
\els 

Due to the length restriction, we omit the proof of this lemma.

\section{A Small-Gain theorem for Time-Delay Systems in State-Space
  Form}\label{s-3} 

Consider a large-scale interconnected system composed of $n$
subsystems:
\beq
& &
\dot x_1(t) = f_1\left((x_1)_t, \  ( v_2)_t, \  ( v_3)_t,  \ldots, ( v_k)_t, \ 
u_1(t)\right),\quad\nonumber\\
& &\dot x_2(t) = f_2\left((x_2)_t, \  ( v_1)_t, \  ( v_3)_t,
\ldots, ( v_k)_t, \  u_2(t)\right),\nonumber\\
& &\phantom{\dot x_2(t)} \vdots\label{e-sys-int}\\
& &\dot x_k(t) = f_k\left((x_k)_t, \  (v_1)_t, \ ( v_2)_t, \ldots, 
(v_{k-1})_t,  u_k(t)\right),\nonumber
\eeq
subject to the interconnection
\be{e-interc}
v_i = x_i,
\quad 1\le i\le k.
\ee
For each $1\le i\le k$, assume
that each $f_i$ is locally Lipschitz jointly on all of its entries,
and for each $i$ and each $t$, 
$
x_i(t) \in\R^{n_i}, \ v_i(t) \in\R^{n_i}$ and
$u_i(t)\in\R^{m_i}$.

For any $\kk$-function $\rho$, we say that $\rho < \id$ if $\rho(s) < s$
for all $s> 0$.

\bt{sg-delay}
Suppose that for each $x_i$-subsystem of \rref{e-sys-int}, there exist
$\kk$-functions $\sigma_i, \gamma_{ij}$ and $\gamma_i^u$ such that the following
properties hold:
\begin{itemize}
\item the GS property: 
\beq
& &\abs{x_i(t)} \leq \max_{i \neq j}  
\left\{
\sigma_i\left(\norm{x_i}_{[-\theta,0]}\right), \right. \nonumber \\ 
& &\quad\left.
\gamma_{ij}\left(\norm {v_j}_{[-\theta,\infty)}\right), 
\gamma_i^u\left(\norm{u}_{[0,\infty)}\right) \right\};
\eeq
for all $ t\ge 0$, and
\item the AG property:
\beq 
& &  \overline{\lim_{t\rightarrow \infty}} \abs {x_i(t)} \leq \max _{i \neq
    j}\left\{\gamma_{ij}\left(\norm {v_j}_{[-\theta,\infty)}\right) \right., \nonumber \\
& & \left. \gamma_i^u\left(\norm{u}_{[0,\infty)}\right) \right\}.
\label{e-tagi}
\eeq
 \end{itemize}
Assume further that the set of cyclic small-gain conditions hold: 
 \be{e-sg}
   \gamma_{i_1i_2}\circ\gamma_{i_2i_3}\circ\cdots\circ\gamma_{i_ri_1} <\id 
 \ee
for all $2\le r\le k$,
$1\le i_j\le k$, $i_j\not=i_{j'}$ if $j\not=j'$.
Then the interconnected system \rref{e-sys-int}-\rref{e-interc} is
forward complete, and it admits
the AG and GS properties with $u=(u_1, \ldots u_k)$ as inputs.
That is, 
there exist class $\kk$-functions ${\sigma}(\cdot),\ 
{\gamma}^u(\cdot)$ such that:
\be{e-tgs2}
\abs{x(t)} \leq \max \left\{\sigma\left(\norm{x}_{[-\theta, 0]}\right),\
\gamma^u\left(\norm{u}_{[0,\infty)}\right)\right\} 
\ee
for all $t\ge 0$, and 
\be{e-tag}
\overline{\lim_{t\rightarrow \infty}} \abs {x(t)}  \leq {\gamma^u}
\left(\norm{u}_{[0,\infty)}\right).
\ee
\ets
\medskip

Note that for any $\rho_1, \rho_2\in\kk$, $\rho_1\circ\rho_2<\id$ if and only
if $\rho_2\circ\rho_1 < \id$.  Consequently, to verify the set of small-gain
conditions \rref{e-sg} for all choices of
$\gamma_{i_1i_2}\circ\gamma_{i_2i_3}\circ\cdots\circ\gamma_{i_ri_1}$ for
which $r\ge 2$, and $1\le i_j\le k$, $i_j\not=i_{j'}$ if $j\not=j'$, it
is sufficient to verify \rref{e-sg} for all choices of those
$\gamma_{i_1i_2}\circ\gamma_{i_2i_3}\circ\cdots\circ\gamma_{i_ri_1}$ 
with $i_1 <  \min\{i_2, \ldots, i_r\}$. 

When $k=2$, the set of small-gain conditions \rref{e-sg} becomes the
usual small-gain condition: $\gamma_{12}\circ\gamma_{21} < \id$.  For
the case of $k=3$, the set of small-gain conditions 
\rref{e-sg} becomes the following:
\beqn
&\gamma_{12}\circ\gamma_{21} < \id,\; \gamma_{13}\circ\gamma_{31}
< \id,\;
\gamma_{23}\circ\gamma_{32} < \id;&\\
&\gamma_{12}\circ\gamma_{23}\circ\gamma_{31} < \id,\;
\gamma_{13}\circ\gamma_{32}\circ\gamma_{21} < \id.&
\eeqn 

In the special case when the subsystems in \rref{e-sys-int} are free of
the external signals $u_i(\cdot)$, the interconnected system becomes
\beq
& &
\dot x_1(t) = f_1\left((x_1)_t, \  ( v_2)_t, \  ( v_3)_t,  \ldots, (
  v_k)_t\right),\quad\nonumber\\ 
& &\phantom{\dot x_2(t)} \vdots\label{e-sys-into}\\
& &\dot x_k(t) = f_k\left((x_k)_t, \  (v_1)_t, \ ( v_2)_t, \ldots, 
(v_{k-1})_t\right),\nonumber
\eeq
subject to the interconnection
\be{e-interco}
v_i = x_i,
\quad 1\le i\le k.
\ee
The  following is then an immediate consequence of Theorem~\ref{sg-delay}:
\bc{c-t1}
Suppose that for each $x_i$-subsystem of \rref{e-sys-into}, there exist
$\kk$-functions $\sigma_i$ and $\gamma_{ij}$ ($j\not= i$, $1\le j\le k$)
such that the following properties hold:
\begin{itemize}
\item the GS property: 
\[
\abs {x_i(t)} \leq \max_{i \neq j}\left\{
\sigma_i\left(\norm{x_i}_{[-\theta,0]}\right),\ 
\gamma_{ij}\left(\norm
  {v_j}_{[-\theta,\infty)}\right)\right\}; 
\]
for all $t\ge 0$, and
\item the AG property:
 \[
   \overline{\lim_{t\rightarrow \infty}} \abs {x_i(t)} \leq \max _{i \neq
    j}\left\{\gamma_{ij}\left(\norm
   {v_j}_{[-\theta,\infty)}\right)\right\}. 
\]
 \end{itemize}
Assume further that the set of cyclic small-gain conditions \rref{e-sg}
holds for all $2\le r\le k$,
$1\le i_j\le k$, $i_j\not=i_{j'}$ if $j\not=j'$.
Then, the interconnected system \rref{e-sys-into}-\rref{e-interco} is
globally asymptotically stable in the following sense:
\begin{itemize}
\item for some $\sigma\in\kk$, $\abs{x(t)}\le \sigma\left(
\norm{x}_{[-\theta, 0]}\right)$ for all $t\ge 0$;
\item $\displaystyle{\lim_{t\rightarrow \infty} \abs {x(t)} =
    0}$.~\mybox 
\end{itemize}
\ecs

\subsection{An Example}
In what follows, we consider an example of a system composed of
three subsystems with delays (without $u$ for simplicity).
Let $\Delta > 0$ be a constant time-delay. 

Consider the system described by the equations,
\beq
\dot{x}_1 (t) &=& -3x_1(t)+ \frac{v_2^2(t-\Delta)}{1+v_2^2(t-\Delta)}, \nonumber \\
\dot{x}_2 (t) &=& \frac{-3}{2}x_2(t)+ v_3^3(t-\Delta), \nonumber \\
\dot{x}_3 (t) &=& -2x_3(t)+ v_1^2(t-\Delta),
\label{ex1} 
\eeq 
with the interconnection $v_i = x_i$ for $i=1, 2, 3$.
Each of the $x_i$-subsystems satisfies the AG and GS conditions. 
The
gain functions can be chosen as:
\comment{\beqn
\gamma_{12} &=& \frac{s^2}{2(1+s^2)},\\
\gamma_{23} &=& s^3,\quad \gamma_{31} = s^2,
\eeqn}
\[
\sigma_1(s) =7s ,\ \sigma_2(s) = 4s,\ \sigma_3 = 3s, \text{and}
\]
\[
\gamma_{12} = \frac{s^2}{2(1+s^2)},\ \gamma_{23} = s^3,\ \gamma_{31} = s^2.
\]
We show the computations for the calculation of $\sigma_1(s)$
and $\gamma_{12}(s)$. The other gain functions can be calculated in a similar manner.
Let $w_1(t)=\frac{v_2^2(t-\Delta)}{1+v_2^2(t-\Delta)}$. 
We can now rewrite the first equation of the system as
\[
\dot{x}_1 (t) = -3x_1(t)+ w_1(t).
\]
The solution of this linear system 
satisfies
\[
\abs{x_1(t)} \leq  \abs{x_{10}}e^{-3t}+\frac{1}{3}\norm{w_1}_{[-\theta,\infty)}.
\]
Using the fact that $a+b < \max\left\{(1+\ve^{-1})a,\ (1+\ve)b\right\}$ for any $\ve>0$ we get that
\[
\abs{x_1(t)} \leq \max \left\{7\abs{x_{10}}e^{-3t},\ \frac{1}{2}\norm{w_1}_{[-\theta,\infty)}\right\}  
\]
by letting $ \ve=\frac{1}{6}$.
To verify the small-gain condition, it is enough to show
that $\gamma_{12}\circ \gamma_{23}\circ \gamma_{31} < \id$.  By
calculation, 
\[
\gamma_{12}\circ \gamma_{23}\circ \gamma_{31}(s)=\frac{s^{12}}{2(1+s^{12})}.
\]
The desired small-gain condition $\frac{s^{12}}{2(1+s^{12})} < s$ is
equivalent to $\frac{s^{12}}{2}<s+s^{13}$. The inequality can be
verified by considering the two cases of $s\le 1$ and $s > 1$.
By the small-gain theorem, the interconnected system given in \rref{ex1}
is globally asymptotically stable.~\qed

\section{Input/Output Operators}

To prove the small-gain theorem for systems of differential equations
with delays as stated in the previous section, we first consider
the more general case of the small-gain theorem for input/output
operators, an approach used in the work \cite{is-smallgain}. Our results
established for operators allow small-gain theorems to be developed
for several situations; systems of differential equations with delays
being just one particular application.

\subsection{Small-Gain Theorem for Operators}

We say that a triple $(\tau, y, u)$ is a trajectory if 
$\tau \in [0,\infty]$, $u=(u_1,...,u_q): [0, \tau)\rightarrow \R^q$ is
measurable and locally
essentially  bounded, and $y=(y_1,...,y_p): [0, \tau)\rightarrow\R^p$ is
continuous.

Note that the trajectories are defined in an abstract way, and
no underlying relation is presupposed between the functions $u$ and $y$.

For an initialized system as in \rref{e-sys} with $x_0(s)=\xi(s)$ on
$[-\theta, 0]$, let $\hat u = (v, u)$, $y(t) = x(t, \xi, v, u)$,
then for any $0<\tau < T^{\max}_{\xi,v,u}$,  the triple
$(\tau, y, \hat u)$ can be identified as a trajectory defined in this section.

\bp{op-prop} Consider a trajectory $(\tau, y, u)$ for which $\tau=~\infty$. 
Assume the following conditions hold:
\begin{itemize}
\item there exist some
class $\kk$-functions $\gamma_{ij}(\cdot)$, $\gamma_i^u(\cdot)$, and
a constant $c\ge 0$, such that
 \be{e-1}
   \abs {y_i(t)} \leq \max_{i \neq j} \left\{c,\ 
    \gamma_{ij}\left(\norm {y_j}_{[0,t)}\right),\ \gamma_i^u\left(\norm{u}_{[0,t)}\right)
  \right\}   
\ee
for all $t\ge 0$;  and
\item it holds that 
\beq
\overline{\lim_{t \rightarrow \infty}} \abs {y_i(t)} & \leq & \max_{i \neq
    j} \left\{\gamma_{ij}\left(\overline{\lim_{t \rightarrow \infty}} \abs {y_j(t)}\right), \right. \nonumber \\
 & &\qquad\quad\left.    
    \gamma_i^u\left(\norm{u}_{[0,\infty)}\right) \right\}. \label{ag}
 \eeq \
\end{itemize} 
Assume further, the small-gain condition: 
 \be{e-smallgain}
   \gamma_{i_1i_2}\circ\gamma_{i_2i_3}\circ\cdots\circ\gamma_{i_ri_1} <\id ,\ 
 \ee
for all $2\le r\le p$, $i_j\not=i_{j'}$ whenever $j\not=j'$.
Then there exist $\kk$-functions $\ \widetilde{\sigma}_i(\cdot),\ 
  \widetilde{\gamma}^u_i(\cdot) $ and $\widehat{\gamma}^u_i
(\cdot)$ such that:
\be{e-p-gs}
\abs{y_i(t)}  \leq \max \left\{\widetilde{\sigma}_i(c),\
  \widetilde{\gamma}^u_i\left(\norm u_{[0,\infty)}\right)\right\} 
\ee
and
\be{e-p-ag}
\overline{\lim_{t \rightarrow \infty}} 
\abs {y_i(t)}  \leq \widehat{\gamma}^u_i
\left(\norm{u}_{[0,\infty)}\right).
\ee
\eps

To prove the proposition we first consider the following lemma.  Note
that this lemma holds for any $\tau \in [0,\infty]$ though the
proposition only requires the lemma for the special case when $\tau
=\infty$. We include in our consideration the case when $\tau < \infty$
to develop a result applicable to interconnected systems as in
Section~\ref{s-3} when forward completeness is not known a priori.

\bl{op1} Consider a trajectory $(\tau, y, u)$, where $\tau \in [0,
\infty]$.  Suppose that for this trajectory, condition \rref{e-1} holds 
for all $t\in [0, \tau)$.  Assume the gain 
functions $\gamma_{ij}$ satisfy the
small-gain condition \rref{e-smallgain}.  Then there exist
$\widetilde{\sigma}_i, \widetilde{\gamma}^u_i \in\kk$ such
that 
\be{e-l1}
\abs{y_i(t)}  \leq \max \left\{\widetilde{\sigma}_i(c),\
  \widetilde{\gamma}^u_i\left(\norm u_{[0,\tau)}\right)\right\},
\ee
for all $t\in [0, \tau)$.
\el

The lemma can be proved by induction on $p$. The case of $p=2$ is in
fact part of the known small-gain theorem for systems with two subsystems
(though not stated in the form for operators). 
 Due to the length restriction, we omit the proof of this lemma.

\medskip

\subsection{Proof of Proposition \ref{op-prop}}

The first part of the proposition follows immediately from Lemma
\ref{op1} with $\tau=\infty$. The proof of the limit property of the proposition is also
inductive.  For the sake of saving space, we skip the proof for the case of $p=2$. 

Instead of treating the general inductive step to pass from $p$
to $p+1$, we just go through the case of $p=3$.  

For each $i, j =1, 2, 3$, let $\gamma_{ij}$ and $\gamma_i^u$ be
$\kk$-functions such that the small-gain condition holds for
$\{\gamma_{ij}\}$. 
Let $(\infty, y, u)$ be a
trajectory satisfying all assumptions of the proposition with the given
gain functions $\{\gamma_{ij}\}$ and $\{\gamma_i^u\}$.  Furthermore,
assume that $\norm{u}_{[0, \infty)}< \infty$.  Let
$b_i$ denote $\displaystyle{\overline{\lim_{t\rightarrow \infty}} \abs
{y_i(t)}}$.  Property \rref{e-p-gs} implies that $b_i< \infty$ for
$i=1, 2, 3$. By the assumption in \rref{ag}, we have:
\[
b_1 \leq \max\left\{\gamma_{12}(b_2), \gamma_{13}(b_3),
  \gamma^u_1\left(\norm{u}_{[0,\infty)}\right)\right\}, 
\]
\[
b_2 \leq \max\left\{\gamma_{21}(b_1), \gamma_{23}(b_3),
  \gamma^u_2\left(\norm{u}_{[0,\infty)}\right)\right\} 
\]
and
\[
b_3 \leq \max\left\{\gamma_{31}(b_1), \gamma_{32}(b_2), 
\gamma^u_3\left(\norm{u}_{[0,\infty)}\right)\right\}.
\]

We eliminate $b_3$ from the first two inequalities as in the
following: 
\beqn
b_1 \leq &\max& \left\{
\gamma_{12}(b_2), \gamma_{13}\circ\gamma_{31}(b_1),
\gamma_{13}\circ \gamma_{32}(b_2), \right.\\
& &\quad\left.
\gamma_{13}\circ\gamma^u_3\left(\norm{u}_{[0, \infty)}\right), 
\gamma^u_1\left(\norm{u}_{[0,\infty)}\right)
\right\}.
\eeqn
\beqn
b_2 \leq &\max& \left\{
\gamma_{21}(b_1), \gamma_{23}\circ \gamma_{31}(b_1),
\gamma_{23}\circ\gamma_{32}(b_2), \right.\\
& & \quad \left.
\gamma_{23}\circ\gamma^u_3\left(\norm{u}_{[0,\infty)}\right), 
\gamma^u_2\left(\norm{u}_{[0,\infty)}\right)
\right\}.
\eeqn
By the small-gain condition, $\gamma_{ij}\circ \gamma_{ji}(s)<s$ 
for all $s > 0$, we get:
\beqn
b_1 \leq &\max& \left\{
\gamma_{12}(b_2),\gamma_{13}\circ \gamma_{32}(b_2),\right.\\
& &\quad\left. 
\gamma_{13}\circ\gamma^u_3\left(\norm{u}_{[0,\infty)}\right), 
\gamma^u_1\left(\norm{u}_{[0,\infty)}\right)\right\},
\eeqn
\beqn
b_2 \leq & \max & \left\{
\gamma_{21}(b_1), \gamma_{23}\circ \gamma_{31}(b_1),\right.\\
& &\quad\left. 
\gamma_{23}\circ\gamma^u_3\left(\norm{u}_{[0,\infty)}\right), 
\gamma^u_2\left(\norm{u}_{[0,\infty)}\right)\right\}.
\eeqn
Now define
\[
\widetilde{\gamma_{ij}}(s):=\max\left\{
\gamma_{ij}(s),\gamma_{i3}\circ \gamma_{3j}(s) \right\},
\]
\[
\widetilde{\gamma_{i}}^u(s):=\max\left\{
\gamma_{i}(s),
\gamma_{i3}\circ \gamma^u_{3}(s) \right\}.
\]
We then have 
\[
b_1 \leq \max\left\{
\widetilde{\gamma_{12}}(b_2), 
\widetilde{\gamma_1}^u\left(\norm{u}_{[0,\ \infty)}\right)\right\},
\]
\[
b_2 \leq \max\left\{
\widetilde{\gamma_{21}}(b_1), 
\widetilde{\gamma_2}^u\left(\norm{u}_{[0,\ \infty)}\right)\right\},
\]
and consequently,
\beqn
b_1 \leq &\max & \left\{
\widetilde{\gamma_{12}}\circ\widetilde{\gamma_{21}}(b_1), 
\widetilde{\gamma_{12}}\circ\widetilde{\gamma_2}^u\left(\norm{u}_{[0,
    \infty)}\right), \right.\\
& &\qquad\quad\left.
\widetilde{\gamma_1}^u\left(\norm{u}_{[0,\infty)}\right)\right\}.
\eeqn
It can be shown that 
$\widetilde{\gamma_{12}}$ and $\widetilde{\gamma_{21}}$ satisfy the
small-gain condition $\widetilde{\gamma_{12}}\circ
\widetilde{\gamma_{21}} < \id$, and thus, 
\[
b_1 \leq \widehat{\gamma}^u_1\left(\norm{u}_{[0, \infty)}\right),
\]
where $\widehat{\gamma}^u_1 (s) = \max\{
\widetilde{\gamma_{12}}\circ\widetilde{\gamma_2}^u(s), 
\widetilde{\gamma_1}^u(s)\}$. Similarly, one sees that 
\[
b_2\le \widehat{\gamma}^u_2\left(\norm{u}_{[0,\infty)}\right),
\]
where $\widehat{\gamma}^u_2(s) = \max\{
\widetilde{\gamma_{21}}\circ\widetilde{\gamma_1}^u(s), 
\widetilde{\gamma_2}^u(s)\}$. 

With the obtained estimates on $b_1$ and $b_2$, one has 
\beqn
b_3\leq &\max &\left\{
\gamma_{31}\left(\widehat{\gamma}^u_1\left(\norm u_{[0,\infty)}\right)\right), \
 \right. \\
& & \quad \left.
\gamma_{32}\left(\widehat{\gamma}^u_2\left(\norm u_{[0,\infty)}\right)\right), \
\gamma^u_3\left(\norm{u}_{[0,\infty)}\right) \right\}.
\eeqn
{}From this it follows that
\[
b_3\leq\widehat{\gamma}^u_3\left(\norm {u}_{[0,\infty)}\right), 
\]
where $\widehat{\gamma}^u_3(s)=\max\left\{
\gamma_{31}\circ\widehat{\gamma}^u_1(s),
\gamma_{32}\circ\widehat{\gamma}^u_2(s),
\gamma^u_3(s) \right\}$. This concludes the proof for the case of
$p=3$.

\subsection{Proof of \protect{Theorem
    \ref{sg-delay}}}\label{s-proof-thm} 

Let $u$ and $\xi = (\xi_1, \ldots, \xi_k)$ be given.  
Consider the corresponding trajectory
$x(t) = (x_1(t), \ldots, x_k(t))$ of the interconnected system
\rref{e-sys-int}-\rref{e-interc} defined on the maximum interval 
$[0,\ T^{\max}_{\xi, u})$. Let $T = T^{\max}_{\xi,u}$.  
Then one has the following on $[0, T)$ for each $i$:
\beqn
\abs {x_i(t)} &\leq& \max_{i \neq j} \left\{
\sigma_i\left(\norm{\xi_i}_{[-\theta,\ 0]}\right),\right.\\ 
& &\qquad\left.
\gamma_{ij}\left(\norm {x_j}_{[-\theta, \ t)}\right), 
\gamma_i^u\left(\norm{u}_{[0,\ T)}\right) \right\}.
\eeqn
Observe that 
\[
\norm{x_j}_{[-\theta,\ t)}\leq \max\left\{\norm{x_j}_{[-\theta,\ 0]},\
  \norm{x_j}_{[0,\ t)}\right\}. 
\]
Hence, 
\beqn
\abs {x_i(t)} &\leq & \max_{i \neq j} \left\{
\sigma_i\left(\norm{\xi_i}_{[-\theta,\ 0]}\right),
 \gamma_{ij}\left(\norm {x_j}_{[-\theta, \ 0]}\right), \right.\\
& &\qquad\left.
\gamma_{ij}\left(\norm {x_j}_{[0, \ t)}\right), \ 
\gamma_i^u\left(\norm{u}_{[0,T)}\right) \right\}.
\eeqn
Let $\displaystyle{c=\max_{i \neq j}\left\{ \sigma_i\left(\norm{\xi_i}_{[-\theta,\ 0]}\right),\ 
    \gamma_{ij}\left(\norm {\xi_j}_{[-\theta, \ 0]}\right)\right\}}$. Then we have
    \[
    \abs {x_i(t)} \leq \max_{i \neq j}\left\{c, \ \gamma_{ij}\left(\norm {x_j}_{[0, \ t)}\right), \ 
    \gamma_i^u\left(\norm{u}_{[0,\ T)}\right) \right\}.
    \] 
And thus we can apply Lemma \ref{op1} to $(T, x, u)$ to get
\be{e-tpgs}
\abs{x_i(t)}  \leq \max \left\{\widetilde{\sigma}_i(c),\
  \widetilde{\gamma}^u_i\left(\norm u_{[0,T)}\right)\right\} 
\ee
for all $t\in[0, T)$ and all $i$. This shows that $x(\cdot)$ remains
bounded on the maximum interval $[0, T)$ if $u$ is essentially
bounded on $[0, T)$. 
 From this it follows that $T=\infty$.  (Otherwise, 
$T < \infty$. Then $u$ is essentially bounded on $[0, T)$, and hence
$x(\cdot)$ remains bounded on $[0, T)$, contradicting the
maximality of $T$.)  This in turn implies that \rref{e-tpgs} holds for
all $0 \le t < T=\infty$. (Note that we have assumed $f$ is completely continuous.
 Related proofs will be included 
in a more detailed version of this work.)

Now we can apply Proposition \ref{op-prop} to the system to get the (AG)
property \rref{e-tag}.
Using Lemma \ref{aglem}, one sees that (AG) condition for the
$x_i$-system is equivalent to
\[
   \overline{\lim_{t\rightarrow \infty}} \abs {x_i(t)} \leq \max _{i \neq
    j}\left\{
\gamma_{ij}\left(\overline{\lim_{t\rightarrow\infty}}\abs{v(t)}\right),
\gamma_i^u\left(\overline{\lim_{t\rightarrow\infty}}\abs{u(t)}\right)
\right\}
\]
Hence, for the composed system \rref{e-sys-int}-\rref{e-interc}, the
following holds:
\[
\overline{\lim_{t\rightarrow \infty}} \abs {x_i(t)} \leq 
\max _{i \neq     j}\left\{ 
\gamma_{ij}\left(\overline{\lim_{t\rightarrow \infty}}
\abs {x_j(t)}\right), 
\gamma_i^u\left(\overline{\lim_{t\rightarrow \infty}}\abs{u(t)}\right)
\right\}.
\]
Now we can apply Proposition \ref{op-prop} to conclude that there exists
some $\widehat{\gamma}_i^u$ such that
\[
\overline{\lim_{t \rightarrow \infty}} \abs {x_i(t)}  \leq 
\widehat{\gamma}^u_i
\left(\norm{u}_{[0,\ \infty)}\right)
\]
for all $i$.~\qed

\section{A Remark on ISS for Time-Delay Systems}\label{s-iss}

Notice that our main result Theorem \ref{sg-delay} was presented in terms
of gain functions in the context of the GS and AG properties instead of
the gain functions appearing in an ISS estimate of the type
\rref{isseq}. In the delay-free case, it is well-known that the ISS
property defined by \rref{isseq} is equivalent to the combination of the GS
and AG properties defined by \rref{gseq} and \rref{ageq} (see
\cite{iss-new}).  With this equivalence relation, one sees that
a small-gain theorem in terms of the gain functions in the GS and AG
properties also leads to a small-gain theorem in terms of ISS gain
functions as in \rref{isseq}.  However, it is not clear at this stage if the
combination of the GS and the AG conditions is equivalent to the ISS
property for time-delay systems.

Consider systems of the following type:
\be{e-sysv}
\dot x(t) = f(x_t, u),\  x_0(\cdot) = \xi(\cdot),
\ee
where $f: \xc\times\R^m\rightarrow\R^n$ is locally Lipschitz.  Assume
further that for any bounded sets $K_x\subset\xc$ and $K_u\subset\R^m$,
$\{f(\xi, u):  \xi\in K_x, u\in K_u\}$ is bounded. (Note that this is
not a trivial property since $K_x$ does not need to be compact.)  Although
it is still unclear as to whether the combination of the
GS and the AG conditions is equivalent to the ISS 
property for time-delay systems, we have nevertheless obtained the
following preparatory results (the proofs of which will be given 
in a more detailed version of this work). 

Consider now a system of the following type:
\be{sigmad}
(\Sigma_d): \ \dot x(t) = F(x_t, d(t)), 
\ee
where the disturbance function
$d: \R_{\ge 0}\rightarrow [0, 1]^m$ is measurable. We assume that $f$ is
a locally Lipschitz map.
Let $\Omega$ be the set of measurable functions $d$ with $\abs{d(t)}\le
1$ for all $t\ge0$.  For each $\xi \in \xc $ and $d\in\Omega$, we let $x(t, \xi,
d)$ denote the trajectory of the system corresponding to the initial
state $\xi$ and the disturbance function $d(\cdot)$.
   
\bd{GAS} 
The system $(\Sigma_d)$ is said to be {\it globally asymptotically
stable} (GAS) if 
\begin{itemize} 
\item[(a)] there exists a $\ki$-function $\sigma$ such that 
\[
\abs{x(t,\xi, d)}\leq \sigma\left(\norm{\xi}_{[-\theta,\ 0]}\right)
\]
for all $t\ge 0$; and
\item[(b)] for each trajectory, it holds that
$\displaystyle{ \lim_{t\rightarrow \infty}\abs{x(t, \xi, d)}=0}$. 
\end{itemize} 
\eds

\bd{UGAS}
 A system as in \rref{sigmad} is {\it uniformly globally asymptotically
   stable} (UGAS) if it satisfies property (a) in Definition \ref{GAS},
 and the following holds:
\[
\forall \veps>0 \ \forall \kappa>0 \ \exists \ 
T=T(\veps,\ \kappa)\geq 0 \ s.t.: 
\]
\[
\norm \xi \leq \kappa\Rightarrow\sup_{t\geq T}\abs{x(t,\xi,d)}\leq \veps
\ \ \forall\,d  \in  \Omega. 
\]
\eds

\addtolength{\textheight}{-0.516in}

Clearly, a system \rref{sigmad} is globally asymptotically stable if it
is uniformly globally asymptotically stable.
 
Let $\varphi:C[-\theta,\ 0] \rightarrow \R_{\geq 0}$ be any locally
Lipschitz functional such that $\varphi({\bf{0}})=0$ where {\bf{0}}
denotes the zero function. Consider the auxiliary system associated with
the system $\displaystyle{(\Sigma_u)}$;  
\be{e-sysf}
(\Sigma_\varphi):\ \dot{x}(t)= f\left(x_t,\ \varphi(x_t)d(t)\right),
\ee
where $d\in\Omega$.

Let $x_{\varphi}(t, \xi, d)$ denote the trajectory of
$(\Sigma_\varphi)$ with initial state $\xi$ and input $d$.

\bp{ag-gs-gas}
Suppose that a system $(\Sigma_u)$ as in \rref{e-sys}
satisfies both the  (AG) and the (GS) properties. 
Then there exists a class $\ki$ function $\rho$ which is locally
Lipschitz such that with $\varphi(\xi) =\rho\left(\norm{\xi}_{[-\theta,
    0]}\right)$, the corresponding auxiliary system $(\Sigma_\varphi)$
as in \rref{e-sysf} is  
globally asymptotically stable.
\ep

\bp{UGASthenISS}  Consider a system $(\Sigma_u)$ as in \rref{e-sys}.
Suppose:
\begin{itemize}
\item the system satisfies the (GS) property; and
\item there is a  class $\ki$-function $\rho$ which is locally Lipschitz
such that with $\varphi(\xi) = \rho\left(\norm{\xi}_{[-\theta,0]}\right)$,
  the corresponding auxiliary
system $\displaystyle{(\Sigma_{\varphi})}$ as in \rref{e-sysf} is
  uniformly globally asymptotically stable.
\end{itemize}
Then the system $\displaystyle{(\Sigma_u)}$ satisfies the ISS property. 
\ep

In order to show that the combination of the (GS) and the (AG)
conditions is equivalent to the ISS property, a crucial step is to
determine for the system  ($\Sigma_\varphi$) if the
global asymptotic stability property implies uniform
global asymptotic stability. This remains a topic
for further study.
 
\section{Conclusion}

The nonlinear ISS small-gain theorem has been generalized to large-scale
systems with time-delays.  Both state-space form and input-output
operators are considered for large-scale system modeling.  Under the set
of cyclic small-gain conditions, it is shown that the large-scale system
enjoys the same type of stability properties as each individual
subsystem.  Our future work will be directed at applications of this
tool to the control of time-delay nonlinear systems. 

\section{Appendix A}
In this section we prove the following:

\bl{l-Lip} Suppose $f:\xc\times\vz\times\R^n\rightarrow\R^n$ is locally
Lipschitz, then for any $v\in\vc$, any $u\in\uc$, the map
\[
F:\xc\times [0, \infty)\rightarrow\R^n, \quad
F(\xi, t) := f(\xi, v_t, u(t))
\]
is locally Lipschitz in $\xi$, uniformly in $t$ for all $t$ in any
compact set. 
\el

{\it Proof.} Let $v\in\vc, u\in\uc$ be given.  Consider a compact
subset $K$ of $\xc$ and a compact interval $[0, T]$.  Since $u$ is
locally essentially bounded, there exists some $L >0$ such that
$\norm{u}_{[0, T]}\le L$. Since $v$ is
continuous, $v$ is uniformly continuous on $[-\theta, T]$.  Thus, for
any $\ve>0$ given, there exists some $\delta>0$ such that for  all $t\in
[0, T]$ 
\[
\abs{v_t(s_1) - v_t(s_2)} 
= \abs{v(t+s_1) - v(t+s_2)} < \ve
\]
for all $s_1, s_2 \in [-\theta, 0]$.  This shows that the family of functions
$\{v_t\}_{0\le t\le T}$ is equicontinuous.  It is clear that the family
$\{v_t\}_{0\le t\le T}$ is uniformly bounded. It can also be shown 
that the family is closed. 
Thus by the Arzel{\`a}-Ascoli theorem, the set ${\cal W}= \{v_t\}_{0\le t\le T}$ is a
compact subset of $\vz$.  Since $f$ is Lipschitz on $K\times{\cal
  W}\times [-L, L]^m$, there exists some $M\ge 0$ such that
\beqn
& &\abs{f(\xi_1, w_1, \mu_1) - f(\xi_2, w_2, \mu_2)} \\
& & \le M(\norm{\xi_1-\xi_2}_{[-\theta, 0]} + \norm{w_1-w_2}_{[-\theta, 0]} + \abs{\mu_1 - \mu_2})
\eeqn
for all $\xi_1, \xi_2 \in K$, all $w_1, w_2 \in{\cal W}$, and all
$\mu_1,\mu_2 \in[-L, L]^m$. In particular, for all  $\xi_1, \xi_2\in K$,
almost all $t\in [0, T]$,
\[
\abs{f(\xi_1, v_t, u(t) ) - f(\xi_2, v_t, u(t))}\le M\norm{\xi_1 -
 \xi_2}_{[-\theta, 0]},
\]
which means 
$\abs{F(\xi_1, t) - F(\xi_2, t)}\le M\norm{\xi_1 - \xi_2}_{[-\theta, 0]}$ for almost
all $t\in [0, T]$.~\qed

\end{document}